\documentclass[a4paper,12pt]{article}

 \usepackage[utf8]{inputenc}
 \usepackage[dvips]{graphicx}
 \usepackage{color}
 \usepackage{graphicx}
 \usepackage{amsmath,amsfonts,amsthm,amssymb}
 \usepackage[english]{babel}
\usepackage{ifpdf}
\usepackage{anysize}
\marginsize{2cm}{2cm}{2cm}{2cm}

\theoremstyle{plain}
\newtheorem{thm}{Theorem}
\newtheorem{lemma}[thm]{Lemma}

\theoremstyle{definition}

\newtheorem{definition*}{Definition}

\theoremstyle{remark}

\renewcommand{\int}{\mathop{\rm int}}
\renewcommand{\epsilon}{\varepsilon}

\newcounter{ris}
\renewcommand{\r}{\refstepcounter{ris}%
                  Fig.~\arabic{ris}.}

\newcounter{upr}

\begin{document}
	
\ifpdf
\DeclareGraphicsExtensions{.pdf, .jpg, .tif, .mps}
\else
\DeclareGraphicsExtensions{.eps, .jpg, .mps}
\fi

\title{The Lemniscate of Bernoulli, without Formulas}
\author{Arseniy V. Akopyan}
\date{}

\maketitle{}
\abstract{In this paper, we give purely geometrical proofs of the well-known properties of the lemniscate of Bernoulli.}

\subparagraph{What is the lemniscate?}

\emph{A polynomial lemniscate with foci $F_1$, $F_2$, \dots, $F_n$} is a locus of points $X$ such that the product of distances from $X$ to the foci is constant ($\displaystyle\prod_{i=1,...,n}|F_iX|=const$).
The $n$-th root of this value is called the \emph{radius} of the lemniscate.
It is clear, that lemniscate is an algebraic curve of degree (at most) $2n$.
You can see a family of lemniscates with three foci in Figure~\ref{fig:Petrunin3points}.

\begin{center}
 \includegraphics[angle=180]{lemniscatePetrunin}\\
\r \label{fig:Petrunin3points} The lemniscate with three foci.
\end{center}


A lemniscate with two foci is called a \emph{Cassini oval}. It is named after the astronomer Giovanni Domenico Cassini who studied them in 1680.
The most well-known Cassini oval is the \emph{lemniscate of Bernoulli} which was described by Jakob Bernoulli in 1694.
This is a curve such that for each point of the curve, the product of distances to foci equals quarter of square of the distance between the foci (Fig.~\ref{fig:lemBernulli}). Bernoulli considered it as a modification of an ellipse, which has the similar definition: the locus of points with the sum of distances to foci is constant. (Bernoulli was not familiar with the work of Cassini).

It is clear that the lemniscate of Bernoulli passes through the midpoint between the foci.
This point is called the \emph{juncture} or \emph{double point} of the lemniscate.

\begin{center}
 \includegraphics[scale=1]{lemniscate-1}
 
\r \label{fig:lemBernulli}   $|F_1X|\cdot |F_2X|=|F_1O|\cdot |F_2O|$.
\end{center}
The lemniscate of Bernoulli has many very interesting properties.
For example, the area bounded by the lemniscate is equal to $\frac {1}{2}|F_1F_2|^2$.
In this paper we prove some other properties, mainly using purely
synthetic arguments.

\subparagraph{How to construct the lemniscate of Bernoulli?}

There exists a very simple method to construct the lemniscate of Bernoulli using the following three-bar linkage.
This construction was invented by James Watt:
Take two equal rods $F_1A$ and $F_2B$ of the length $\frac{1}{\sqrt 2}|F_1F_2|$ and fixed at the points $F_1$ and $F_2$ respectively.
Let points $A$ and $B$ lie on opposite sides of the line $F_1F_2$.
The third rod cinnects the points $A$ and $B$ and its
length equals $|F_1F_2|$ (Fig.~\ref{fig:construct of lemniscate}).
Then, during the motion of the linkage \emph{the midpoint $X$ of the rod $AB$ traces the lemniscate of Bernoulli with foci at $F_1$ and $F_2$}.

\begin{center}
 \includegraphics[scale=1]{lemniscate-2}\\
\r \label{fig:construct of lemniscate}
\end{center}

\begin{proof}
	Note that the quadrilateral $F_1AF_2B$ is an isosceles trapezoid (Fig.~\ref{fig:proof of method}).
	Moreover, triangles $\triangle AF_1X$ and $\triangle ABF_1$ are similar, because they have the common angle $A$ and the following relation on their sides holds:
	\[
	\frac{|AF_1|}{|AX|}=\frac{|AB|}{|AF_1|}=\sqrt 2.
	\]

	For the same reason, triangles  $\triangle BXF_2$ and $\triangle BF_2A$ are similar. They have the common angle $B$, and the ratios of the length of the sides with endpoints at $B$ is $\sqrt 2$.
	Therefore, we can write the following equations:
	 \[
	 \angle AF_1X = \angle ABF_1= \angle BAF_2= \angle XF_2B\text{.}
	 \]
	\begin{center}
	 \includegraphics[scale=1]{lemniscate-3}\\
	\r \label{fig:proof of method}
	\end{center}

	Let us remark that in the trapezoid $F_1AF_2B$ angles $\angle A$ and $\angle F_2$ are equal.
	Since angles $\angle XAF_2$ and $\angle XF_2B$ are equal too, we obtain $\angle F_1AX=\angle XF_2A$.
	This implies that triangles $\triangle F_1AX$ and $\triangle AF_2X$ are similar.
	Therefore,
	\[
	\frac{|F_1X|}{|AX|}=\frac{|AX|}{|XF_2|} \Rightarrow |XF_1|\cdot |XF_2|=|AX|^2= |F_1O|^2\text{.}
	\]
	Thus, we have shown that point $X$ lies on the lemniscate of Bernoulli

	Since the motion of the point $X$ is continuous and $X$ 
	attains the farthest points of the lemniscate, the trajectory of $X$ is the whole lemniscate of Bernoulli.
\end{proof}


\begin{center}
 \includegraphics[scale=1]{lemniscate-4}\\
\r \label{fig:prove of maclaurin}
\end{center}

Let $O$ be the midpoint of the segment $F_1F_2$ (double point of the lemniscate).
Denote by $M$ and $N$ the midpoints of the segments $F_1A$ and $F_1B$
respectively (Fig.~\ref{fig:prove of maclaurin}).
Translate the point $O$ by the vector $\overrightarrow{NF_1}$.
Denote the new point by $O'$.
Observe  that triangles $\triangle F_1MO'$ and $\triangle NXO$ are congruent.
Moreover, the following equation holds:
\[
|F_1M|=|F_1O'|=\frac{1}{\sqrt2} |F_1O|\text{.}
\]
In other words, the points $M$ and $O'$ lie on the
circle with center at $F_1$ and radius $\frac{1}{\sqrt2}|F_1O|$.

Using this observation, we can obtain another elegant method to construct the lemniscate of Bernoulli.
\begin{center}
 \includegraphics[scale=1]{lemniscate-5}\\
\r \label{fig:maclaurin}
\end{center}

Let us construct the circle with center at one of the foci and radius  $\frac{1}{\sqrt2}|F_1O|$.
On each secant $OAB$ (where $A$ and $B$ are the points of intersection of the circle and the secant) chose points $X$ and $X'$ such that $|AB|=|OX|=|OX'|$ (Fig.~\ref{fig:maclaurin}).
The union of all points $X$ and $X'$ form the lemniscate of Bernoulli with foci $F_1$ and~$F_2$.

\begin{center}
 \includegraphics[scale=1]{lemniscate-6}\\
\r \label{fig:another way}
\end{center}

Another interesting way to construct lemniscate with the linkages given in Figure~\ref{fig:another way}.
The lengths of the segments $F_1A$ and $F_1O$ are equal.
The point $A$ is the vertex of rods $AX$ and $AY$ of the length $\sqrt{2}|F_1O|$.
Denote the midpoints of these rods by $B$ and $C$, and join them with $O$ by the another rod of the length $\frac{|AX|}{2}$.
In the process of rotating of point $A$ around the circle each of the points $X$ and $Y$ generates a half of the lemniscate of Bernoulli with foci $F_1$ and $F_2$.


%
%

\subparagraph{The lemniscate of Bernoulli and equilateral hyperbola.}

The hyperbola is a much more well-known curve.
Hyperbola with foci $F_1$ and $F_2$ is a set of all points $X$ such that value  $\big| |F_1X|-|F_2X|\big|$ is constant.
Points $F_1$ and $F_2$ are called the foci of the hyperbola.
Among all hyperbolas, we single out \emph{equilateral} or \emph{rectangular hyperbolas}, i. e.,  the set of points $X$ such that $\big||F_1X|-|F_2X|\big|=\frac{|F_1F_2|}{\sqrt{2}}$.


The lemniscate of Bernoulli is an inversion image of an equilateral hyperbola.
Before proving this claim, let us recall the definition of a inversion.

\begin{definition*}
	Inversion with respect to the circle with center $O$ and radius $r$ is a transformation which maps every point $X$ in the plane to the point $X^*$ lying on the ray $OX$ such that $|OX^*|=\frac{r^2}{|OX|}$.
\end{definition*}

The inversion has many interesting properties, see, for example~\cite{johnson2007advanced}.
Among the properties, is the following: a circle or a line will invert to either a circle or a line, depending on whether it passes through the origin.
We will prove here just one simple lemma that will help us later.

\begin{lemma}
\label{inversion of the line}
Suppose $A$ is an orthogonal projection of the point $O$ on some line~$\ell$.
Then inversion image of the line $\ell$ with respect to a circle with center at $O$ is the circle with diameter $OA^*$, where $A^*$ is the inversion image of the point~$A$.
\end{lemma}

\begin{center}
 \includegraphics[scale=1]{lemniscate-7}\\
\r \label{fig:inversion}
\end{center}

\begin{proof}
 Let $B$ be any point on the line $\ell$, and let $B^*$ be its image (Fig.~\ref{fig:inversion}).
 Since
 \[
|OA^*|=\frac{r^2}{|OA|} \text{ \ \ and \ \ }
|OB^*|=\frac{r^2}{|OB|}\text{,}
 \]
 we obtain that triangles $\triangle OAB$ and $\triangle OB^*A^*$ are similar.
Therefore the angle $\angle OB^*A^*$ is a right angle and the point $B^*$ lies on the circle with diameter $OA^*$.
\end{proof}

Note that \emph{the center $O^*_1$ of this circle is the inversion of the point $O_1$, where $O_1$  is the point symmetric to $O$ with respect to the line $\ell$}.

Now, let us prove that \emph{the lemniscate of Bernoulli with foci $F_1$ and $F_2$, is an inversion of the equilateral hyperbola with foci $F_1$ and $F_2$ with respect to the circle with center at $O$ and radius $|OF_1|$}.
For this proof, we will use the results we obtained in the proof of correctness of the first method to construct the lemniscate (Fig.~\ref{fig:proof of method}).
Let $P$ be the point of the intersection of the lines $F_1A$ and $F_2B$ and
let $Q$ be symmetric point to $P$ with respect to the line $F_1F_2$ (Fig.~\ref{fig:hyperandber}).
Note that
\begin{align*}
 |F_2Q|-|F_1Q|=|F_2P|-|F_1P|=|AP|-|F_1P|=|F_1A|=\frac{|F_1F_2|}{\sqrt{2}}\text{.}
\end{align*}

Therefore, the points $P$ and $Q$ lie on the equilateral hyperbola with foci at $F_1$ and~$F_2$.
Now it remains to show that points $X$ and $Q$ are the images of each other under the inversion with center at $O$ and radius $|OF_1|$.
First, let us show that triangles $\triangle F_1XO$ and $\triangle PF_1O$ are similar.

\begin{center}
 \includegraphics[scale=1]{lemniscate-8}\\
\r \label{fig:hyperandber}
\end{center}

The quadrilateral $F_1XOB$ is a trapezoid. Therefore $\angle OXF_1+\angle XF_1B=180^\circ$.
Also, we have $\angle AF_1O+\angle OF_1P=180^\circ$.
Since the angle $\angle XF_1B$ is equal to the angle $\angle AF_1O$, we obtain that  $\angle OXF_1$ and $\angle OF_1P$ are equal to each other.

Since angels $\angle XF_2B$ and $\angle XF_1A$ are equal, we have that $\angle XF_1P+\angle PF_2X=180^\circ$.
In other words, the quadrilateral $PF_1XF_2$ is inscribed.
Therefore, we have
\[
\angle F_2F_1X=\angle F_2PX=\angle F_1PO\text{.}
\]
The last equation provided that points $O$ and $X$ are symmetric to each other with respect to the perpendicular bisector of the segment $F_1B$.

Thus, triangles $\triangle F_1XO$ and $\triangle PF_1O$ are similar because they have two corresponding pairs of equal angles.
It follows that $\angle F_1OX$ and $\angle F_1OP$ are equal, and we have that the point $Q$ lies on the ray $OX$.
In addition, from similarity of triangles $\triangle F_1XO$ and $\triangle QF_1O$ (it is congruent to the triangle $\triangle PF_1O$), we obtain
\[
\frac{|OX|}{|OF_1|}=\frac{|OF_1|}{|OQ|}
\Rightarrow
|OX|\cdot |OQ|=|OF_1|^2\text{.}
\]
It means that points $Q$ and $X$ are images of each other under the inversion with center at $O$ and radius $|OF_1|$.

\begin{center}
 \includegraphics[scale=1]{lemniscate-9}\\
\r \label{fig:inversion_of_tangent}
\end{center}

If we look at Figure~\ref{fig:hyperandber}, we can make another observation: the points $X$ and $O$ lie on the circle with centered at~$P$.
It is interesting that \emph{this circle touches
the lemniscate of Bernoulli.}

\begin{proof}
Suppose $\ell$ is a tangent line to the hyperbola in the point $Q$.
From Lemma \ref{inversion of the line}, it follows that the image of the line $\ell$ under the inversion with center at $O$ and radius~$|F_1O|$ is a circle $\omega_\ell$ passing through the point $O$.
Since~$X$ is the inversion image of the point $Q$, we see that the circle $\omega_\ell$ touches the lemniscate at the point~$X$.
From the same Lemma we conclude that the center of this circle lies on the normal line from the point $O$ to the line $\ell$.
Let us show that lines $OP$ and $OQ$ are symmetric to each other with respect to the line $F_1F_2$.
It follows that the point~$P$ is the center of the circle $\omega_\ell$.

\begin{center}
 \includegraphics[scale=1]{lemniscate-10}\\
\r \label{fig:symmetry}
\end{center}

Without loss of generality, we can assume that the equation of the hyperbola is $y=\frac{1}{x}$.
Suppose line $\ell$ intersects the abscissa and the ordinate in the points $R$ and $S$ respectively (Fig.~\ref{fig:symmetry}).
It is well-known that the derivative of the function $\frac{1}{x}$ at the point $x_0$ is equal to $\frac{-1}{x_0^2}$.
It follows that the point $Q$ is the midpoint of the segment $RS$, and $OQ$ is the median of the right triangle $\triangle ROS$. Therefore, the angles $\angle QOR$ and $\angle QRO$ are equal.
Since angles $\angle POS$ and $\angle QOR$ are also equal, we obtain that the lines~$OP$ and $RS$ are perpendicular, as was to be proved.
\end{proof}

\begin{center}
 \includegraphics[scale=1]{lemniscate-11}\\
\r \label{fig:normal}
\end{center}

Let us note the following: Since the circle $\omega_\ell$ touches the lemniscate at the point~$X$, the radius $PX$ of this circle is a normal (perpendicular to the tangent line) to the lemniscate at $X$ (Fig.~\ref{fig:normal}).
Note that the triangle $\triangle XPO$ is isosceles, and the lines $XO$ and $PO$ are symmetric with respect to the line $F_1O$.
Therefore, we can write these equations:
\[
\angle PXO=\angle XOP=2\angle POF_1\text{.}
\]
It follows that there exists the following very simple method to construct the normal to the lemniscate of Bernoulli:
For any point $X$ on the lemniscate, take the line forming with the line intersecting $OX$ at $X$, with the angle of intersection equal to $2\angle XOF$.
This line will be a normal to the lemniscate.


\end{document}